\documentclass[a4paper,10pt]{article}
\usepackage{kubik}
\hypersetup{pdftitle={Curvature bounded below: a definition a la Berg--Nikolaev}%
,pdfauthor={Nina Lebedeva and Anton Petrunin}
}
\begin{document}

\title{Curvature bounded below:\\
a definition a la Berg--Nikolaev}
\author{N.~Lebedeva and A.~Petrunin}
\date{}
\maketitle

\section{Introduction}

In this note we give a new characterization of spaces with curvature $\ge 0$ in the sense of Alexandrov.
Our work is inspired by \cite{berg-nikolaev} and \cite{sato}, where an analogous definition was given for curvature $\le 0$.

\begin{thm}{Main theorem}
Let $\mathcal X$ be a complete space with intrinsic metric.
Then $\mathcal X$ is an Alexandrov space with curvature $\ge 0$ 
if and only if any quadruple $p,x,y,z\in \mathcal X$ 
satisfies the following inequality
$$|p x|^2+|p y|^2+|p z|^2
\ge
\tfrac13\cdot\!\l(|x y|^2+|y z|^2+|z x|^2\r).
\eqlbl{*}$$

\end{thm}

The inequality \ref{*} is  quite weak.
For example, one can%
\footnote{Say, consider the  metric on $\{p,x,y,z\}$ defined as $|p x|=p y|=|p z|=1$, $|x y|=|x z|=2$ and $|y z|=\eps$ where $\eps>0$ is sufficiently small;
it satisfies \ref{*} for each relabeling but not \ref{1+3}.} 
construct a metric space $\mathcal F$ with 4 points 
which satisfies \ref{*} for each relabeling by $p,x,y,z$, 
such that $\mathcal F$ does not admit an isometric embedding into 
any Alexandrov space with curvature $\ge 0$.

The similar conditions which were known before simply describe all 4-point sets in a nonnegatively curved space.
For instance the following inequality for model angles:
$$\angk{}{p}{x}{y}+\angk{}{p}{y}{z}+\angk{}{p}{z}{x}\le2\cdot\pi.
\eqlbl{1+3}$$
In fact, if a 4-point metric space satisfies \ref{1+3} for each relabeling
then it can be isometrically embedded into 
Euclidean plane 
or a $2$-sphere of some radius $R>0$ (the proof is left to the reader).

\parbf{Why do we care.}
Since the condition \ref{*} is so weak, 
it should be useful as a test to check that a given space has curvature $\ge0$ in the sense of Alexandrov.
However, we are not aware of a single case when it makes life easier.

To explain the real reason why we are interested in this topic, 
we need to reformulate our Main theorem using the language similar to one given in \cite[Section 1.19$_+$]{gromov}.

Let us denote by $\mathbf{M}^4$ the set of isometry classes of 4-point metric spaces.
Let $\mathfrak A$ and $\mathfrak B$ be the sets of isometry classes in $\mathbf{M}^4$ 
which satisfy correspondingly \ref{*} and \ref{1+3} for all relabeling of points by $p,x,y,z$.
(As it mentioned above, $\mathfrak B\subsetneq\mathfrak A$.)
Further, given a metric space $\mathcal X$,
denote by $\mathbf{M}^4(\mathcal X)$ the set of isometry classes of 4-point subspaces in $\mathcal X$.

Main theorem says that if the space $\mathcal X$ has intrinsic metric and $\mathbf{M}^4(\mathcal X)\subset \mathfrak A$ then $\mathbf{M}^4(\mathcal X)\subset \mathfrak B$.
From above, the set $\mathfrak B$ is the smallest set which satisfies the above property for any $\mathcal X$.

It would be interesting to find a general pattern of such phenomena.
Assume you start with arbitrary $\mathfrak A\subset \mathbf{M}^4$, can you figure out what is the corresponding $\mathfrak B$,
or can one describe the properties of $\mathfrak B$ which might appear this way?

Note that Globalization theorem (see \cite{BGP})
as well as Berg--Nikolaev characterization of $\mathrm{CAT}(0)$ spaces 
both admit  interpretations in the above terms.
Also, in \cite{FOERTSCH-LYTCHAK-SCHROEDER}, it was shown that set defined by Ptolemy inequality can appear as $\mathfrak B$.

\parbf{Acknowledgment.}
We want to thank Alexander Lytchak for an interesting discussion in the M\"unster's botanical garden.

\section{The proof}

The ``only if'' part follows  directly from generalized Kirszbraun's theorem, see \cite{lang-schroeder}.
One only has to check that the inequality \ref{*} holds in the model plane.
(Alternatively, one can prove \ref{1+3}~$\Rightarrow$~\ref{*} directly.)

\parit{``if'' part.}
We may assume that $\mathcal X$ is geodesic, otherwise pass to its ultraproduct.

It is sufficient to show that if $z$ is a midpoint of geodesic $[p q]$ in $\mathcal X$ then
$$2\cdot|x z|^2
\ge
|x p|^2+|x q|^2-\tfrac12\cdot|p q|^2,\eqlbl{**}$$
for any $x\in \mathcal X$.

Directly from \ref{*} we have the following weaker estimate%
$$
3\cdot|x z|^2\ge |x p|^2+|x q|^2-\tfrac12\cdot |p q|^2
\eqlbl{***}
$$

\begin{wrapfigure}{r}{22mm}
\begin{lpic}[t(-6mm),b(0mm),r(0mm),l(0mm)]{pics/mediana(0.5)}
\lbl[r]{0,0;$p$}
\lbl[l]{42,1;$q$}
\lbl[l]{34,47;$x_0$}
\lbl[l w]{26,16;$x_1$}
\lbl[l w]{24,6;$x_2$}
\lbl[t]{21,-1;$z$}
\end{lpic}
\end{wrapfigure}

Set $x_0=x$ and consider a sequence of points
 $x_0,x_1,\dots$ on $[x z]$ such that $|x_n z|=\tfrac1{3^n}\cdot|x z|$.
Let $\alpha_n$ be a sequence of real numbers such that
$$\alpha_n\cdot|x_n z|^2
=
|x_n p|^2+|x_n q|^2-\tfrac12\cdot|p q|^2.$$
Applying \ref{*}, we get
$$|x_{n+1}p|^2+|x_{n+1}q|^2+|x_{n+1}x_{n}|^2
\ge
\tfrac13\cdot(|x_n p|^2+|x_n q|^2+|p q|^2).$$
Subtract $\tfrac12\cdot|p q|^2$ from both sides;
after simplification you get
$$\alpha_{n+1}\ge 3\cdot\alpha_n-4.$$

Now assume \ref{**} does not hold; i.e., $\alpha_0>2$
then $\alpha_n\to\infty$ as $n\to\infty$.
On the other hand, from \ref{***}, we get $\alpha_n\le 3$, a contradiction.
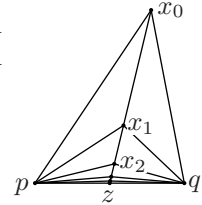
\qeds

\section*{P.S.: Arbitrary curvature bound}

One can obtain analogous characterization of Alexandrov spaces with curvature $\ge\kappa$ for any $\kappa\in\RR$.

Here are the inequalities for cases $\kappa=1$ and  $-1$ which correspond to \ref{*} for quadruple $p,x^1,x^2,x^3$:
$$\l(\sum_{i=1}^3\cos|p x^i|\r)^2
\le
\sum_{i,j=1}^3\cos|x^ix^j|=
3+2\cdot\sum_{i<j}\cos|x^ix^j|.
\eqno\text{\ref{*}}^{+}$$
$$\l(\sum_{i=1}^3\cosh|p x^i|\r)^2
\ge
\sum_{i,j=1}^3\cosh|x^ix^j|
=3+2\cdot\sum_{i<j}\cosh|x^ix^j|;
\eqno\text{\ref{*}}^{-}
$$
Note that in both cases we have equality if $p$ is the point of intersections of medians of  the triangle $[x^1x^2x^3]$ in the corresponding model plane.
(In the model planes, medians also pass through incenter.)

The proof goes along the same lines.
The case $\kappa=1$ also follows from Main theorem 
applied to the cone over the space.

\end{document}